\documentclass[12pt]{article}

\usepackage{amsmath, amssymb}
\usepackage{graphicx}
\usepackage{hyperref}

\newtheorem{theorem}{Theorem}[section]
\newtheorem{corollary}{Corollary}[section]

\numberwithin{equation}{section}
\usepackage{mathrsfs}

\title{Hineva Inequality on Some Submanifolds of Quaternionic Space forms}
\author{Idrees Fayaz Harry, Mehraj Ahmad Lone, Lokenath Ganguly \\ 
	\texttt{harryidrees.96@nic.in, mehraj.jmi@gail.com, lokenathganguly@gmail.com}}
\date{}

\begin{document}
	
	\maketitle
	
	\begin{abstract}
			In this article, we establish Hineva inequality for different types of submanifolds of Quaternionic Space forms.
	\end{abstract}

\section{Introduction} 

The seminal work of John F. Nash \cite{NASH} on the isometric immersion of Riemannian manifolds into Euclidean spaces provides a powerful motivation to view Riemannian manifolds as submanifolds of Euclidean spaces. Consequently, establishing sharp, simple relationships between the intrinsic and extrinsic invariants of a Riemannian submanifold has become a central problem in submanifold theory. Among extrinsic invariants, the squared mean curvature is particularly significant. Intrinsic invariants, on the other hand, include classical quantities such as sectional curvature, Ricci curvature, and scalar curvature.

Introduced by Prof. B. Y. Chen in his seminal 1993 paper \cite{1}, the theory of $\delta$-invariants (or Chen invariants) has grown into one of the most active and important research areas in submanifold differential geometry. The original motivation for these new Riemannian invariants was to address an open problem posed by S.S. Chern concerning the existence of minimal isometric immersions of Riemannian manifolds into Euclidean spaces of arbitrary dimension \cite{2}. Before Chen's work, this problem remained unsolved because classical intrinsic invariants were insufficient for controlling the extrinsic geometry of submanifolds (see \cite{3}).

In 1984, Hineva established lower and upper bounds for the sectional curvatures of submanifolds of a Riemannian manifold in terms of the mean curvature vector and the squared norm of the second fundamental form (see \cite{hineva1985submanifolds}, \cite{hineva1994note}). Subsequently, in 1990, she presented (without proof) similar lower and upper bounds for the Ricci curvature of submanifolds \cite{hineva1990inequalities}.

\begin{theorem}[A] \label{thm:A} \cite{hineva1990inequalities}
	Let $(M,g)$ be an $n$-dimensional submanifold of a Riemannian manifold $(\tilde{M},\tilde{g})$. Then
	\begin{align}
		K(\Pi_p^2) &\geq \tilde{K}(\Pi_p^2) + \frac{n^2|H|^2}{2(n-1)} - \frac{1}{2} \|\sigma\|^2, \\
		K(\Pi_p^2) &\leq \tilde{K}(\Pi_p^2) + \frac{4-n}{2} \|H\|^2 + \frac{n-2}{2n} \|\sigma\|^2 + \sqrt{\frac{2(n-2)}{n}} \|H\|^2 (\|\sigma\|^2 - n\|H\|^2), \\
		\operatorname{Ric}(X) &\leq \tilde{\operatorname{Ric}}_{T_p M}(X) + \frac{n^2|H|^2}{4}, \\
		\operatorname{Ric}(X) &\geq \tilde{\operatorname{Ric}}_{T_p M}(X) + 2(n-1) \|H\|^2 - \frac{n-1}{n} \|\sigma\|^2 &\nonumber\\-& \frac{n(n-2)}{n-1} \sqrt{\frac{(n-1)}{n}} \|H\|^2 (\|\sigma\|^2 - n\|H\|^2)&.
	\end{align}
\end{theorem}

In 2008, she provided a proof for the Ricci curvature bounds $(A_3)$ and $(A_4)$ and discussed the equality cases \cite{hineva2008submanifolds}.\\

In \cite{1}, Chen derived a necessary condition for the existence of minimal isometric immersions of a Riemannian manifold $M$ into Euclidean space. He also established a sharp inequality involving intrinsic invariants (scalar and sectional curvatures) and the principal extrinsic invariant (the squared mean curvature). Later, in \cite{4}, Chen obtained further sharp inequalities relating the $k$-Ricci curvature, the squared mean curvature, and the shape operator for submanifolds of real space forms of arbitrary codimension and analyzed the equality cases.
\begin{theorem}[B] \cite{4} 
	Let $M$ be an $n$-dimensional submanifold of a real space form $\tilde{M}(c)$. Then the following statements are true.
\end{theorem}

\noindent \textbf{(a)} For \( X \in T_p^1 M \), we have
\[
\operatorname{Ric}(X) \leq (n-1)c + \frac{n^2}{4} \|H\|^2.
\]

\textbf{(b)} If \( H(p) = 0 \), then the equality case in \textbf{(B)} is true for \( X \in T_p^1 M \) if and only if \( X \in \mathcal{N}(p) \).

\textbf{(c)} The equality case in \textbf{(B)} is true for every \( X \in T_p^1 M \) if and only if either \( p \) is a geodesic point or \( n = 2 \) and \( p \) is an umbilical point.\\
Notably, these inequalities are optimal, as they are attained by many important classes of submanifolds. Since then, an extensive literature has developed around Chen invariants and related inequalities for various submanifold classes in diverse ambient geometries, including complex space forms \cite{5}--\cite{7}, cosymplectic space forms \cite{8}--\cite{10}, Sasakian space forms \cite{11}--\cite{13}, locally conformal Kähler space forms [\cite{14}--\cite{16}], generalized complex space forms [\cite{17}--\cite{19}], locally conformal almost cosymplectic manifolds \cite{20,21}, $(k,\mu)$-contact space forms \cite{22, 23}, Kenmotsu space forms \cite{24, 25}, S-space forms \cite{26, 27}, and T-space forms \cite{28}; see also \cite{29} and the references therein.

In 2005, using optimization techniques, Oprea obtained a Ricci curvature inequality for submanifolds of real space forms and improved it for Lagrangian submanifolds of complex space forms \cite{oprea2005,OPERA2010}. That same year, Hong and Tripathi extended Chen's work \cite{4} to derive a Ricci curvature inequality for submanifolds of arbitrary Riemannian manifolds, also discussing the equality cases.

\begin{theorem}[C] \label{thm:C} (see \cite{HONGTRIPATHI})
	Let $(M,g)$ be an $n$-dimensional submanifold of a Riemannian manifold $(\tilde{M},\tilde{g})$. Then the following statements are true.
	
	\textbf{(a)} For \( X \in T_p^1 M \), we have
	\[
	\operatorname{Ric}(X) \leq \frac{1}{4} n^2 \|H\|^2 + \tilde{\operatorname{Ric}}_{(T_p M)}(X).
	\]
	
	\textbf{(b)} The equality case in \textbf{(C)} is true for \( X \in T_p^1 M \) if and only if
	\[
	\begin{cases}
		\sigma(X,Y) = 0 & \text{for all } Y \in T_p M \text{ orthogonal to } X, \\
		2\sigma(X,X) = nH(p).
	\end{cases}
	\]
	If \( H(p) = 0 \), then the equality case in \textbf{(C)} is true for \( X \in T_p^1 M \) if and only if \( X \in \mathcal{N}(p) \).
	
	\textbf{(c)} The equality case in \textbf{(C)} is true for all unit vectors \( X \in T_p M \) if and only if either \( p \) is a geodesic point or \( n = 2 \) and \( p \) is an umbilical point.
\end{theorem}

By 2008, inequality \textbf{(C)} was referred to as the Chen--Ricci inequality \cite{tripathiricci}. In 2009, using essential algebraic inequalities, Deng obtained an improved Chen--Ricci inequality for Lagrangian submanifolds of complex space forms and analyzed its equality cases \cite{49}. In \cite{tripathiimproved}, the author derived Chen--Ricci and improved Chen--Ricci inequalities for Kulkarni--Nomizu tensor fields on Riemannian manifolds satisfying the algebraic Gauss equation, also discussing the equality cases. Following Hineva's contributions, inequality \textbf{(A4)} will henceforth be referred to as the Hineva inequality.

\begin{theorem}\cite{hineva2008submanifolds}
	Let \((M,g)\) be an \(n\)-dimensional submanifold of an \(m\)-dimensional Riemannian manifold \((\tilde{M},\tilde{g})\). Then the following statements are true.
	
	\begin{enumerate}
		\item[(a)] For \(X \in T_p^1 M\), we have
		\begin{eqnarray}\label{HINEVAINE}
			\operatorname{Ric}(X) &\geq& \tilde{Ric(T_p M)(X)}\nonumber\\& &+ \frac{n-1}{n} \left( 2n \|H\|^2 - (n-2) \sqrt{\frac{n \|H\|^2 (\|\sigma\|^2 - n \|H\|^2)}{n-1}} \right).
		\end{eqnarray}
		
		\item[(b)] The equality case in \eqref{HINEVAINE} is true for \(X \in T_p^1 M\) if and only if \(p\) is a quasi-umbilical point, that is, the matrices \((\sigma_1^0)\) are of the following form:
		\begin{eqnarray}\label{HINEVAMAT}
			\begin{pmatrix}
				\lambda_n & 0 & \cdots & 0 \\
				0 & \mu_n & \cdots & 0 \\
				\vdots & \vdots & \ddots & \vdots \\
				0 & 0 & \cdots & \mu_n \\
				0 & 0 & \cdots & 0
			\end{pmatrix},
		\end{eqnarray}
		and for the matrices
		\begin{eqnarray}\label{HINEVAMAT1}
			\frac{\lambda_n - \mu_n}{\lambda_n + (n-1)\mu_n}
		\end{eqnarray}
		is invariant over \(\alpha \in \{n+1,\ldots,m\}\), where \(\lambda_n\) and \(\mu_n\) are two different eigenvalues of the matrices \((\sigma_1^0)\) of the second fundamental form and
		\begin{eqnarray}\label{HINEVAMAT2}
			\lambda_n = \frac{1}{n} \left( \operatorname{tr}(\sigma^n) \pm (n-1) \sqrt{\frac{n \operatorname{tr}(\sigma^n)^2 - (\operatorname{tr}(\sigma^n))^2}{n-1}} \right),
		\end{eqnarray}
		\begin{eqnarray}\label{HINEVAMAT3}
			\mu_n = \frac{1}{n} \left( \operatorname{tr}(\sigma^n) \pm \sqrt{\frac{n \operatorname{tr}(\sigma^n)^2 - (\operatorname{tr}(\sigma^n))^2}{n-1}} \right).
		\end{eqnarray}
		
		\item[(c)] The equality case in \eqref{HINEVAINE} is true for every \(X \in T_p^1 M\) if and only if \(p\) is an umbilical point.
	\end{enumerate}
	
\end{theorem}

Recently the Chen--Ricci and Hineva inequalities have been formulated for a Kulkarni--Nomizu tensor field on a Riemannian manifold satisfying the algebraic Gauss equation, with the equality cases investigated \cite{VERMA}. Applying these results, the inequalities are further studied for Riemannian submanifolds. Furthermore, a new formulation for generic submanifolds of almost contact metric manifolds is proposed, and the Hineva inequality is derived for various types of Riemannian submanifolds of generalized Sasakian space forms, including Sasakian, cosymplectic, and Kenmotsu space forms.

\begin{theorem}\cite{VERMA}
	Let \((M, g)\) be an \(n\)-dimensional submanifold of an \(m\)-dimensional Riemannian manifold \((\tilde{M}, \tilde{g})\). Then the following statements are true.\\
	(a)		For \(X \in T_p^1 M\), we have
	
	\begin{eqnarray}\label{tripathi}
		&\text{Ric}_{(T_p M)}(X) + \frac{n-1}{n} \left( 2n \|H\|^2 - \|\sigma\|^2 - (n-2)(n-1) \|H\|^2 \right) \\
		&\leq \text{Ric}(X) \leq \text{Ric}_{(T_p M)}(X) + \frac{n^2}{4} \|H\|^2.
	\end{eqnarray}
	
	(b)	 The simultaneous equality case of both the inequalities in \eqref{tripathi} is true for \(X \in T_p^1 M\) if and only if \(p\) is a quasi-umbilical point that is, the matrices \((h_{ij}^{\alpha})\) are of the following form:
	\begin{eqnarray}\label{tripathimat}
		\begin{pmatrix}
			\lambda_\alpha & 0 & \cdots & 0 & 0 \\
			0 & \mu_\alpha & \cdots & 0 & 0 \\
			\vdots & \vdots & \ddots & \vdots & \vdots \\
			0 & 0 & \cdots & \mu_\alpha & 0 \\
			0 & 0 & \cdots & 0 & \mu_\alpha
		\end{pmatrix},
	\end{eqnarray}
	where \(\lambda_\alpha\) and \(\mu_\alpha\) are two different eigenvalues of each one of the matrices \((h_{ij}^{\alpha})\) such that \(\lambda_\alpha = \operatorname{tr}(h^{\alpha})/2\) and \(\mu_\alpha = \operatorname{tr}(h^{\alpha})/(2(n-1))\).\\
	(c)	The simultaneous equality case of both the inequalities in \eqref{tripathi} is true for every \(X \in T_p^1 M\) if and only if either \(p\) is a geodesic point or \(n=2\) and \(p\) is an umbilical point.

\end{theorem}

\section{Preliminaries}

Let \(\bar{M}^{m}\) be a \(4m\)-dimensional Riemannian manifold with metric \(g\). \(\bar{M}^{m}\) is called a \textit{quaternion Kaehlerian manifold} if there exists a 3-dimensional vector bundle \(V\) of tensors of type \((1,1)\) over \(\bar{M}^{m}\) with local basis of almost Hermitian structures \(I\), \(J\) and \(K\) such that

\begin{enumerate}
	\item[(a)] \(IJ = -JI = K\), \(JK = -KJ = I\), \(KI = -IK = J\), \(I^{2} = J^{2} = K^{2} = -1\).
	\item[(b)] for any local cross-section \(\phi\) of \(V\), \(\bar{\nabla}_{X}\phi\) is also a cross-section of \(V\), where \(X\) is an arbitrary vector field on \(\bar{M}^{m}\) and \(\nabla\) the Riemannian connection on \(\bar{M}^{m}\).
\end{enumerate}

In fact, condition (b) is equivalent to the following condition:

\begin{enumerate}
	\item[(b')] there exist local 1-forms \(p\), \(q\) and \(r\) such that
	\[
	\bar{\nabla}_{X}I = r(X)J - q(X)K, \quad 
	\bar{\nabla}_{X}J = -r(X)I + p(X)K, \quad 
	\bar{\nabla}_{X}K = q(X)I - p(X)J.
	\]
\end{enumerate}		
Now let \(X\) be a unit vector on \(\bar{M}^{m}\), then \(X, IX, JX, KX\) form an orthonormal frame on \(\bar{M}^{m}\). We denote by \(Q(X)\) the 4-plane spanned by them. For any two orthonormal vectors \(X, Y\) on \(\bar{M}^{m}\), if \(Q(X)\) and \(Q(Y)\) are orthogonal, the plane \(\pi(X,Y)\) spanned by \(X, Y\) is called a \textit{totally real plane}. Any 2-plane in a \(Q(X)\) is called a \textit{quaternionic plane}. The sectional curvature of a quaternionic plane \(\pi\) is called the \textit{quaternionic sectional curvature} of \(\pi\). A quaternion Kaehlerian manifold is a \textit{quaternion space form} if its quaternionic sectional curvatures are equal to a constant. A quaternion projective space, denoted by \(QP^{m}(4c)\), is a quaternion Kaehlerian manifold of constant quaternionic sectional curvature \(4c\).\\		
It is known that a quaternionic Kaehlerian manifold \(\bar{M}^{m}\) is a quaternion space form if and only if its curvature tensor \(\bar{R}\) is of the following form \cite{32}:		
\begin{equation}\label{eq:R-curvature}
	\begin{aligned}
		\bar{R}(X,Y)Z &= c\bigl\{g(Y,Z)X - g(X,Z)Y \\
		&\quad + g(IY,Z)IX - g(IX,Z)IY + 2g(X,IY)IZ \\
		&\quad + g(JY,Z)JX - g(JX,Z)JY + 2g(X,JY)JZ \\
		&\quad + g(KY,Z)KX - g(KX,Z)KY + 2g(X,KY)KZ\bigr\},
	\end{aligned}
\end{equation}
for vectors \(X,Y,Z\) tangent to \(\bar{M}^{m}\).\\		
Let \(M^{n}\) be an \(n\)-dimensional Riemannian manifold isometrically immersed in a quaternion projective space \(QP^{m}(4c)\). We denote by \(K(\pi)\) the sectional curvature of \(M\) associated with a plane section \(\pi \subset T_{p}M\), \(p \in M\). Also, we denote by \(h\) the second fundamental form and \(R\) the Riemannian curvature tensor of \(M\). We call \(M^{n}\) a \textit{totally real submanifold} of \(QP^{m}(4c)\) if each 2-plane of \(M^{n}\) is mapped into a totally real plane in \(QP^{m}(4c)\). Consequently, if \(M^{n}\) is a totally real submanifold of \(QP^{m}(4c)\), then \(\phi(TM^{n}) \subset T^{\perp}M^{n}\) for \(\phi = I, J\) or \(K\), where \(T^{\perp}M^{n}\) is the normal bundle of \(M^{n}\) in \(QP^{m}(4c)\).\\		
We know that if \(M^{n}\) is a totally real submanifold of \(QP^{m}(4c)\), then for any orthonormal vectors \(X, Y\) in \(M^{n}\), the plane \(\pi(X,Y)\) spanned by \(X\) and \(Y\) is totally real in \(QP^{m}(4c)\), \(Q(X)\) and \(Q(Y)\) are orthogonal and \(g(X,\phi Y) = g(\phi X,Y) = 0\) for \(\phi = I,J\) or \(K\).\\		
Then the equation of Gauss is given by		
\begin{equation}\label{eq:Gauss}
	\bar{R}(X,Y,Z,W) = R(X,Y,Z,W) + g(h(X,W),h(Y,Z)) - g(h(X,Z),h(Y,W)),
\end{equation}
for any vectors \(X,Y,Z,W\) tangent to \(M\).\\
We denote by \(H\) the mean curvature vector, i.e.,		
\begin{equation}\label{eq:mean-curvature}
	H(p) = \frac{1}{n}\sum_{i = 1}^{n} h(e_i, e_i), \quad p \in M,
\end{equation}
where \(\{e_1,\dots,e_n\}\) is an orthonormal basis of the tangent space \(T_p M\).\\		
Also, we denote by
\[
h_{ij}^r = g(h(e_i, e_j), e_r),
\]
and
\[
\|h\|^2 = \sum_{i,j = 1}^{n} g(h(e_i, e_j), h(e_i, e_j)).
\]
For any \(p \in M\) and for any \(X \in T_p M\), we put
\[
IX = P_1 X + F_1 X, \quad JX = P_2 X + F_2 X, \quad KX = P_3 X + F_3 X,
\]
where \(P_i X \in T_p M\), \(F_i X \in T_p^\perp M\).\\		
Let \(\{e_1,\dots,e_n\}\) be an orthonormal basis of \(T_p M\). We denote by
\[
\|P_l\|^2 = \sum_{i,j = 1}^{n} g^2(P_l e_i, e_j), \quad l = 1,2,3.
\]

\begin{theorem}\cite{43}
	Let $M$ be an $n$-dimensional submanifold in a quaternion projective space $QP^m(4c)$. Then:\\
	(i) For each unit vector $X\in T_pM$, we have
	\begin{eqnarray}
		Ric(X)\leq \frac{n^2}{4}	\|H\|^2+(n-1)c+\frac{3}{2}c\sum_{i=1}\|P_i\|^2
	\end{eqnarray}
\end{theorem}

\section{Main Results On submanifolds In quaternion Projective space $QP^m(4c)$}
The following estimation of	the Ricci curvature for  submanifolds in quaternionic projective space forms was firstly proved in \cite{43}
\begin{theorem}\cite{43}
	Let $M$ be an $n$-dimensional submanifold in a quaternion projective space $QP^m(4c)$. Then:\\
	(i) For each unit vector $X\in T_pM$, we have
	\begin{eqnarray}
		Ric(X)\leq \frac{n^2}{4}	\|H\|^2+(n-1)c+\frac{3}{2}c\sum_{i=1}\|P_i\|^2
	\end{eqnarray}
\end{theorem}
\noindent
Based on the foregoing theroem we obtain the following results

\begin{theorem}
	Let $M$ be an $n$-dimensional submanifold in a quaternion projective space $QP^m(4c)$. Then:\\
	(a) For each unit vector $X\in T_pM$, we have
	\begin{eqnarray}\label{proj1}
		Ric(X)&\geq& (n-1)c+\frac{3}{2}c\sum_{i=1}\|P_i\|^2\nonumber\\&+&\frac{n-1}{n}\left(2n\|H\|^2-\|h\|^2-(n-2)\sqrt{\frac{n\|H\|^2(\|h\|^2-n\|H\|^2)}{n-1}}\right)\nonumber\\
	\end{eqnarray}
	(b) The equality case in \eqref{proj1} is true for $X\in T_p^1M$ if and only if the matrices $h_{ij}^{\alpha}$ takes the form
	\begin{equation}
		\left(
		\begin{array}{cccccc}
			\lambda_{\alpha} & 0 & \cdots & 0 & 0 \\
			0 & \mu_{\alpha} & \cdots & 0 & 0 \\
			\vdots & \vdots & \ddots & \vdots & \vdots \\
			0 & 0 & \cdots & \mu_{\alpha} & 0 \\
			0 & 0 & \cdots & 0 & \mu_{\alpha}
		\end{array}
		\right),
	\end{equation}
	and for the matrices
	\begin{equation}
		\frac{|\lambda_{\alpha}-\mu_{\alpha}|}
		{\lambda_{\alpha}+(n-1)\mu_{\alpha}}
		\quad \text{is invariant over } 
		\alpha \in \{n+1,\ldots,2m\},
	\end{equation}
	where $\lambda_{\alpha}$ and $\mu_{\alpha}$ are two different eigenvalues
	of the matrices $\left(h_{ij}^{\alpha}\right)$ of the second
	fundamental form, and
	\begin{equation}
		\lambda_{\alpha}= \frac{1}{n}\left(\operatorname{tr}(h^{\alpha})
		\mp (n-1)\sqrt{	\frac{n\,\operatorname{tr}(h^{\alpha})^2
				-(\operatorname{tr}(h^{\alpha}))^2}{n-1}}
		\right),
	\end{equation}
	\begin{equation}
		\mu_{\alpha}
		= \frac{1}{n}\left(
		\operatorname{tr}(\sigma^{\alpha})
		\pm \sqrt{
			\frac{n\,\operatorname{tr}(h^{\alpha})^2
				-(\operatorname{tr}(h^{\alpha}))^2}{n-1}}
		\right).
	\end{equation}
	
	(c) The equality case in \eqref{proj1} is true for every
	$X \in T_p^1 M$ if and only if $p$ is an umbilical point.

\end{theorem}
\begin{corollary}
	Let $M$ be an $n$-dimensional submanifold in a quaternion projective space $QP^m(4c)$. Then:\\
	(a) For each unit vector $X\in T_pM$, we have
	\begin{eqnarray}\label{proj2}
		&(n-1)c+\frac{3}{2}c\sum_{i=1}\|P_i\|^2+\frac{n-1}{n}\left(2n\|H\|^2-\|h\|^2-(n-2)\sqrt{\frac{n\|H\|^2(\|h\|^2-n\|H\|^2)}{n-1}}\right)&\nonumber\\
		&\leq	Ric(X)\leq (n-1)c+\frac{3}{2}c\sum_{i=1}\|P_i\|^2+\frac{n^2}{4}	\|H\|^2&
	\end{eqnarray}
	(b) The equality case in \eqref{proj2} is true for $X\in T_p^1M$ if and only if the matrices $h_{ij}^{\alpha}$ takes the form
	\begin{equation}
		\left(
		\begin{array}{cccccc}
			\lambda_{\alpha} & 0 & \cdots & 0 & 0 \\
			0 & \mu_{\alpha} & \cdots & 0 & 0 \\
			\vdots & \vdots & \ddots & \vdots & \vdots \\
			0 & 0 & \cdots & \mu_{\alpha} & 0 \\
			0 & 0 & \cdots & 0 & \mu_{\alpha}
		\end{array}
		\right),
	\end{equation}
	
	where $\lambda_{\alpha}$ and $\mu_{\alpha}$ are two different eigenvalues
	of the matrices $\left(h_{ij}^{\alpha}\right)$ such that $\lambda_{\alpha}=\frac{tr(h^{\alpha})}{2}$ and $\mu_{\alpha}=\frac{tr(h^{\alpha})}{2(n-1)}.$
	
	(c) The equality case in \eqref{proj2} is true for every
	$X \in T_p^1 M$ if and only if $p$ is geodesic point or n=2 and p is an umbilical point.
	
\end{corollary}
\section{Hineva Inequality on CR-Submanifolds of Quaternionic Space Forms}
A submanifold $M$ of a quaternion Kähler manifold $(M,\sigma,g)$ is said to be a \emph{quaternionic CR-submanifold} if there exists two orthogonal complementary distributions $D$ and $D^\perp$ on $M$ such that $D$ is invariant under quaternionic structure and $D^\perp$ is total real (see \cite{44}). An estimation of the Ricci curvature of a quaternionic CR-submanifold in a quaternionic space form has been established in \cite{45}, as follows:

\begin{theorem}\cite{45}
	Let $M$ be an $n$-dimensional quaternion CR-submanifold of a quaternion space form $\tilde{M}(c)$. Then:

	\begin{enumerate}
		\item[(i)] For each unit vector $X \in D_x^{\perp}$, we have
		\begin{equation}\label{CR1}
			\operatorname{Ric}(X) \le (n-1)\frac{c}{4} + \frac{n^2}{4} H^2.
		\end{equation}
		
		\item[(ii)] For each unit vector $X \in D_x$, we have
		\begin{equation}\label{CR2}
			\operatorname{Ric}(X) \le (n+8)\frac{c}{4} + \frac{n^2}{4} H^2.
		\end{equation}
		
		\item[(iii)]
	\end{enumerate}
	
\end{theorem}

\begin{theorem}
	Let $M$ be an $n$-dimensional quaternion CR-submanifold of a quaternion space form $\tilde{M}(c)$. Then:

	\begin{enumerate}
		\item[(i)] For each unit vector $X \in D_x^{\perp}$, we have
		\begin{equation}\label{CRM1}
			\operatorname{Ric}(X) \ge (n-1)\frac{c}{4}+\frac{n-1}{n}\left(2n\|H\|^2-\|h\|^2-(n-2)\sqrt{\frac{n\|H\|^2(\|h\|^2-n\|H\|^2)}{n-1}}\right).
		\end{equation}
		
		\item[(ii)] For each unit vector $X \in D_x$, we have
		\begin{equation}\label{CRM2}
			\operatorname{Ric}(X) \geq (n+8)\frac{c}{4} +\frac{n-1}{n}\left(2n\|H\|^2-\|h\|^2-(n-2)\sqrt{\frac{n\|H\|^2(\|h\|^2-n\|H\|^2)}{n-1}}\right).
		\end{equation}
		\item[(iii)] The equality case in \eqref{CRM1} and \eqref{CRM2} is true for $X\in T_p^1M$ if and only if the matrices $h_{ij}^{\alpha}$ takes the form \eqref{HINEVAMAT}-\eqref{HINEVAMAT3}.
		
		\item[(iv)] If $H(x)=0$, then a unit tangent vector $X$ at $x$ satisfies the equality case of \eqref{CRM1} (respectively \eqref{CRM2}) if and only if
		\[
		X \in D_x^{\perp} \cap N_x
		\quad \text{(respectively } X \in D_x \cap N_x\text{)}.
		\]
	\end{enumerate}
\end{theorem}	

\begin{theorem}
	Let $M$ be an $n$-dimensional quaternion CR-submanifold of a quaternion space form $\tilde{M}(c)$. Then:
	\begin{enumerate}
		\item[(i)] For each unit vector $X \in D_x^{\perp}$, we have
		\begin{eqnarray}\label{CRTRIPATHI1}
			(n-1)\frac{c}{4}+\frac{n-1}{n}\left(2n\|H\|^2-\|h\|^2-(n-2)\sqrt{\frac{n\|H\|^2(\|h\|^2-n\|H\|^2)}{n-1}}\right)\nonumber\\\leq \operatorname{Ric}(X) \le (n-1)\frac{c}{4} + \frac{n^2}{4} H^2.
		\end{eqnarray}
		\item[(ii)] For each unit vector $X \in D_x$, we have
		\begin{eqnarray}\label{CRTRIPATHI2}
			(n+8)\frac{c}{4} +\frac{n-1}{n}\left(2n\|H\|^2-\|h\|^2-(n-2)\sqrt{\frac{n\|H\|^2(\|h\|^2-n\|H\|^2)}{n-1}}\right)\nonumber\\\leq\operatorname{Ric}(X) \le (n+8)\frac{c}{4} + \frac{n^2}{4} H^2.
		\end{eqnarray}
		\item[(iii)] The simultaneous equality case of both the inequalities in \eqref{CRTRIPATHI1} and \eqref{CRTRIPATHI2} is true for \(X \in T_p^1 M\) if and only if \(p\) is a quasi-umbilical point that is, the matrices \((h_{ij}^{\alpha})\) are of the following form \eqref{tripathimat}
		
		\item[(iv)] If $H(x)=0$, then a unit tangent vector $X$ at $x$ satisfies the equality case of \eqref{1.7} (respectively \eqref{1.8}) if and only if
		\[
		X \in D_x^{\perp} \cap N_x
		\quad \text{(respectively } X \in D_x \cap N_x\text{)}.
		\]
	\end{enumerate}
\end{theorem}	
\section{Hineva Inequality on Slant subamifolds of Quaternionic Space Forms}
It is clear that, although quaternionic CR-submanifolds are also the
generalization of both quaternionic and totally real submanifolds, there exists no inclusion
between the two classes of quaternionic CR-submanifolds and slant submanifolds.
The following estimation of the Ricci curvature for slant submanifolds in quaternionic
space forms was firstly proved in \cite{40} , and an alternative nice proof can be found in \cite{47} .
\begin{theorem}\cite{40}
	Let $M$ be a $\theta$-slant proper submanifold of a quaternionic space form $\tilde{M}(c)$. Then:
	
	\item[(i)] For each unit vector $X\in T_pM$, we have
	\begin{eqnarray}
		Ric(X)\leq(n-1)\frac{c}{4}+\frac{3c}{8}cos^2\theta+\frac{n^2}{4}	\|H\|^2
	\end{eqnarray}

	\item[(ii)] If $H(x)=0$, then a unit tangent vector $X$ at $x$ satisfies the equality case of  if and only if $X$ belongs to the relative nulln space of M at p
	\begin{eqnarray}
		N_p=\{X\in T_pM / h(X,Y)=0, \forall Y \in T_pM\}
	\end{eqnarray}
	\item[(iii)] The equality case in the theorem holds identically for all unit tangent vectors at p if and only p either p is a totally geodesic point or n=2 and p is totally umbilical point.
\end{theorem}	
\begin{theorem}
	Let $M$ be a $\theta$-slant proper submanifold of a quaternionic space form $\tilde{M}(c)$. Then:
	
	\item[(i)] For each unit vector $X\in T_pM$, we have
	\begin{eqnarray}\label{slanthineva1}
		&Ric(X)\geq& (n-1)\frac{c}{4}+\frac{3c}{8}cos^2\theta\nonumber\\&&+\frac{n-1}{n}\left(2n\|H\|^2-\|h\|^2-(n-2)\sqrt{\frac{n\|H\|^2(\|h\|^2-n\|H\|^2)}{n-1}}\right)\nonumber\\
	\end{eqnarray}
	\item[(ii)] The equality case in \eqref{slanthineva1} is true for $X\in T_p^1M$ if and only if the matrices $h_{ij}^{\alpha}$ takes the form	\eqref{HINEVAMAT}-\eqref{HINEVAMAT3}
	\item[(ii)] If $H(x)=0$, then a unit tangent vector $X$ at $x$ satisfies the equality case of  if and only if $X$ belongs to the relative nulln space of M at p
	\begin{eqnarray}
		N_p=\{X\in T_pM / h(X,Y)=0, \forall Y \in T_pM\}
	\end{eqnarray}
	\item[(iii)] The equality case in the theorem holds identically for all unit tangent vectors at p if and only p either p is a totally geodesic point or n=2 and p is totally umbilical point.
\end{theorem}	
\begin{theorem}
	Let $M$ be a $\theta$-slant proper submanifold of a quaternionic space form $\tilde{M}(c)$. Then:
	
	\item[(i)] For each unit vector $X\in T_pM$, we have
	\begin{eqnarray}\label{thetaslant2}
		(n-1)\frac{c}{4}+\frac{n-1}{n}\left(2n\|H\|^2-\|h\|^2-(n-2)\sqrt{\frac{n\|H\|^2(\|h\|^2-n\|H\|^2)}{n-1}}\right)\nonumber\\\leq \operatorname{Ric}(X) \le (n-1)\frac{c}{4} + \frac{n^2}{4} H^2+\frac{3c}{8}cos^2\theta
	\end{eqnarray}
	\item[(ii)] The simultaneous equality case of both the inequalities in \eqref{thetaslant2} is true for \(X \in T_p^1 M\) if and only if \(p\) is a quasi-umbilical point that is, the matrices \((h_{ij}^{\alpha})\) are of the following form \eqref{tripathimat}.
	\item[(iii)] The simultaneous equality case in the theorem holds identically for all unit tangent vectors at p if and only p either p is a totally geodesic point or n=2 and p is totally umbilical point.
\end{theorem}

\noindent	
{\bf{Conflict of interests:}} The authors declare that they have no conflict of interest, regarding the publication of this paper.

\noindent
\textbf{Acknowledgement:}\\

\end{document}